\documentclass[a4paper, 11pt, twoside]{article}
\usepackage{mfpaperstuff}

\begin{document}
\newcommand\po{{+}\negthinspace1}
\newcommand\nn[1]{N_{\doms#1}}
\newcommand\ind{\operatorname{Ind}}
\newcommand\res{\operatorname{Res}}
\newcommand\iind{\operatorname{-Ind}}
\newcommand\ires{\operatorname{-Res}}
\newcommand\syd{\,\triangle\,}
\newcommand\mo{-1}
\newcommand\mq{-Q}
\newcommand\reg{\mathcal{R}}
\newcommand\swed[1]{\stackrel{#1}{\wed}}
\newcommand\inv[2]{\llbracket#1,#2\rrbracket}
\newcommand\one{\mathbbm{1}}
\newcommand\bsm{\begin{smallmatrix}}
\newcommand\esm{\end{smallmatrix}}
\newcommand{\rt}[1]{\rotatebox{90}{$#1$}}
\newcommand\id{\operatorname{identity}}
\newcommand{\ol}{\overline}
\newcommand{\ul}{\check}
\newcommand\partn{\mathcal{P}}
\newcommand{\py}[3]{\,_{#1}{#2}_{#3}}
\newcommand{\pyy}[5]{\,_{#1}{#2}_{#3}{#4}_{#5}}
\newcommand{\sss}{\mathfrak{S}_}
\newcommand{\nin}{\notin}
\newcommand{\nchar}{\operatorname{char}}
\newcommand{\thmcite}[2]{\textup{\textbf{\cite[#2]{#1}}}\ }
\newcommand\zez{\mathbb{Z}/e\mathbb{Z}}
\newcommand\zepz{\mathbb{Z}/(e+1)\mathbb{Z}}
\newcommand\zo{\bbn_0}
\newcommand\dw{^\triangle}
\newcommand\wod{^\triangledown}
\newcommand{\hhh}{\mathcal{H}_}
\renewcommand{\bbb}{\mathcal{B}_}
\newcommand{\aaa}{\mathcal{A}_}
\newcommand{\sect}[1]{\section{#1}}
\newcommand{\ff}{\mathfrak{f}}
\newcommand{\fff}{\mathfrak{F}}
\newcommand\cf{\mathcal{F}}
\newcommand\sx{x}
\newcommand\bra[1]{|#1\ran}
\newcommand\arb[1]{\widehat{\bra{#1}}}
\newcommand\foc[1]{\mathcal{F}_{#1}}
\newcommand\wed\wedge
\newcommand\wede\barwedge
\newcommand\uu[1]{\,\begin{array}{|@{\,}c@{\,}|}\hline #1\\\hline\end{array}\,}
\newcommand{\ux}[1]{\operatorname{ht}_{#1}}
\newcommand\mmod{\ \operatorname{Mod}}
\newcommand\cgs\succcurlyeq
\newcommand\cls\preccurlyeq
\newcommand\inc{\mathfrak{A}}
\newcommand\fsl{\mathfrak{sl}}
\newcommand\ba{\mathbf{s}}
\newcommand\ta{\tilde\ba}
\newcommand\kt[1]{|#1\rangle}
\newcommand\tk[1]{\langle#1|}
\newcommand\ket[1]{s_{#1}}
\newcommand\jn\diamond
\newcommand\UU{\mathcal{U}}
\newcommand\MM[1]{M^{\otimes#1}}
\newcommand\add{\operatorname{add}}
\newcommand\rem{\operatorname{rem}}
\newcommand\lra\longrightarrow
\newcommand\tru[1]{{#1}_-}
\newcommand\ste[1]{{#1}_+}
\newcommand\lad{\mathcal{L}}
\newcommand\hsl{\widehat{\mathfrak{sl}}}
\newcommand\GG{H}
\newcommand\dable{restrictable}
\newcommand\infi[1]{$(\infty,#1)$-irreducible}
\newcommand\infr[1]{$(\infty,#1)$-reducible}
\newcommand\infs[1]{$#1$-signature}
\newcommand\be[2]{B^{#1}(#2)}
\newcommand\domi{dominant}

\title{Irreducible Specht modules for\\Iwahori--Hecke algebras of type $B$}
\runninghead{Irreducible Specht modules for Iwahori--Hecke algebras of type $B$}
\toptitle

\begin{abstract}
We consider the problem of classifying irreducible Specht modules for the Iwahori--Hecke algebra of type $B$ with parameters $Q,q$.  We solve this problem completely in the case where $q$ is not a root of unity, and in the case $q=-1$ we reduce the problem to the corresponding problem in type $A$.
\end{abstract}

\section{Introduction}

Let $\bbf$ be a field, and $q$ a non-zero element of $\bbf$.  Let $\aaa n$ denote the Iwahori--Hecke algebra of type $A_{n-1}$ over $\bbf$ with parameter $q$.  This algebra arises in various mathematical contexts, and its representation theory closely resembles the modular representation theory of the symmetric group $\sss n$.  For every partition $\la$ of $n$, there is a module $S^\la$ for $\aaa n$ called the \emph{Specht module}.  Specht modules are important because they arise as cell modules for a particular choice of cellular basis for $\aaa n$.  In the case where $\aaa n$ is semisimple, the Specht modules are precisely the irreducible $\aaa n$-modules; in general, the irreducible $\aaa n$-modules arise as the cosocles of certain Specht modules.  It is an interesting question to ask which Specht modules are irreducible, and this question has been answered in almost all cases in a series of papers \cite{jmjs,jmp2,slred,mfred,mfirred,jlm,slcp}.  In particular, the question of which ordinary irreducible representations of the symmetric group remain irreducible in prime characteristic is solved.  For the remaining case in the classification of irreducible Specht modules, some progress has been made in \cite{fl,mfprime,fl2}.

In the present paper we begin the study of irreducible Specht modules in type $B$.  Let $Q$ be an element of $\bbf$, and let $\bbb n$ denote the Iwahori--Hecke algebra of type $B_n$ over $\bbf$ with parameters $Q,q$.  For this algebra there is a cellular basis and a Specht module theory, but with Specht modules indexed by bipartitions of $n$.  Again, we consider the problem of classifying the irreducible Specht modules.  We solve this problem in the fairly easy case where $q$ is not a root of unity in $\bbf$.  Then we consider the case where $q=\mo$; we show how to reduce the main problem in this case to the (still unsolved) type $A$ problem.

In a future paper we hope to consider the remaining cases, where $q$ is an $e$th root of unity for $e>2$.  We also hope to address a generalisation to higher levels.  Specifically, we hope to consider irreducible Specht modules for the \emph{cyclotomic Hecke algebra} or \emph{Ariki--Koike algebra}, which has the same definition as the Hecke algebra of type $B$, but with the quadratic relation for the generator $T_0$ replaced by an arbitrary polynomial relation. The representation theory of this algebra admits a very similar combinatorial framework to $\bbb n$, but with bipartitions replaced by multipartitions (with the number of components being the degree of the relation for $T_0$).  To each such multipartition is associated a Specht module, and one can ask which Specht modules are irreducible.  Some of the results in this paper will generalise fairly easily; the main obstruction is that we do not have an explicit formula for decomposition numbers in the case where $q$ is not a root of unity.

In the next section, we give a brief account of the background results we need.  In Section \ref{infsec}, we address the case $e=\infty$, and in Section \ref{e2sec} we look at the case $e=2$.

\section{Background and basic results}

In this section we briefly treat some background on the representation theory of Iwahori--Hecke algebras.  An excellent reference for Hecke algebras of type $A$ is the book by Mathas \cite{math}.  For type $B$, the paper of Dipper and James \cite{dj} gives a useful introduction.

\subs{Iwahori--Hecke algebras of types $A$ and $B$}

Throughout this paper, $\bbf$ denotes a field, and $q,Q$ are elements of $\bbf$ with $q\neq0,1$.  For a positive integer $n$, the \emph{Iwahori--Hecke algebra} $\bbb n$ is the unital associative $\bbf$-algebra with generators $T_0,\dots,T_{n-1}$ and relations
\begin{alignat*}2
(T_0-Q)(T_0+1) &=0&&\\
(T_i-q)(T_i+1) &=0&&(1\ls i\ls n-1)\\
T_0T_1T_0T_1&=T_1T_0T_1T_0&&\\
T_iT_{i+1}T_i&=T_{i+1}T_iT_{i+1}&\qquad&(1\ls i\ls n-2)\\
T_iT_j&=T_jT_i&&(0\ls i\ls j-2\ls n-3).
\end{alignat*}
We denote by $\aaa n$ the subalgebra generated by $T_1,\dots,T_{n-1}$; this is the Iwahori--Hecke algebra of type $A_{n-1}$.

We denote by $e\in\{2,3,\dots\}\cup\{\infty\}$ the multiplicative order of $q$ in $\bbf$.

\subs{Partitions and bipartitions}

Recall that a \emph{partition} of $n$ is a weakly decreasing sequence $\la=(\la_1,\la_2,\dots)$ of non-negative integers with $|\la|:=\la_1+\la_2+\dots=n$.  When we write a partition, we usually omit trailing zeroes and group together equal parts with a superscript.  We write $\varnothing$ for the unique partition of $0$.

A \emph{bipartition} of $n$ is an ordered pair $(\la,\mu)$ of partitions with $|\la|+|\mu|=n$.  We refer to $\la$ and $\mu$ as the first and second \emph{components} of $(\la,\mu)$.

\subsubs{Young diagrams}

The \emph{Young diagram} $[\la]$ of a partition $\la$ is the set
\[\lset{(i,j)}{j\ls\la_i}\subset\bbn^2.\]
We refer to elements of $\bbn^2$ as \emph{nodes}, and to elements of $[\la]$ as nodes of $\la$.  A node $\fkn$ of $\la$ is \emph{removable} if $[\la]\setminus\{\fkn\}$ is the Young diagram of a partition, while a node not in $[\la]$ is an \emph{addable node} of $\la$ if $[\la]\cup\{\fkn\}$ is the Young diagram of a partition.  We use the English convention for drawing Young diagrams, in which $i$ increases down the page and $j$ increases from left to right.  For example, the following illustrates the Young diagram of the partition $(4^2,2)$, with removable nodes marked $-$ and addable nodes marked $+$.
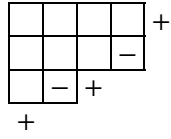
\[\gyoung(;;;;:+,;;;;-,;;-:+,:+)\]

The Young diagram $[\la,\mu]$ of a bipartition $(\la,\mu)$ is the subset
\[\lset{(i,j)_1}{j\ls\la_i}\cup\lset{(i,j)_2}{j\ls\mu_i}\]
of $\bbn^2\times\{1,2\}$; we use the terms node, addable node and removable node as for partitions.

Given $q,Q$ as in the definition of $\bbb n$, we define the \emph{residue} of a
node $(i,j)_k$ to be
\[
\begin{cases}
-Qq^{j-i}&(k=1)\\
\phantom{-Q}q^{j-i}&(k=2).
\end{cases}
\]
For $x\in\bbf$, a node of residue $x$ will be called an
\emph{$x$-node}.

\subsubs{The dominance order}

The \emph{dominance order} on partitions is a partial order $\dom$ defined by saying that $\la\dom\xi$ if for every $i\gs1$ we have $\la_1+\dots+\la_i\gs\xi_1+\dots+\xi_i$.  The dominance order on bipartitions (also denoted $\dom$) is defined by $(\la,\mu)\dom(\xi,\nu)$ if $\la\dom\xi$ and $|\la|+\mu_1+\dots+\mu_i\gs|\xi|+\nu_1+\dots+\nu_i$ for every $i\gs1$.

\subsubs{Beta-sets}

A \emph{beta-set} (also referred to as a \emph{one-runner abacus}) is a subset $B$ of $\bbz$ such that for $N\gg0$ we have $-N\in B\notni N$.  Given a partition $\la$ and an integer $i$, one may define a beta-set $B^i(\la)$ by
\[B^i(\la)=\lset{\la_j+i-j}{j\gs1}.\]
Conversely, any beta-set $B$ defines a unique partition: writing the elements of $B$ as $b_1>b_2>\dots$, there is a unique integer $i$ such that $b_N=i-N$ for $N\gg0$; then $B=B^i(\la)$, where $\la_j=b_j+j-i$ for each $j$.

It is conventional to depict a beta-set with an abacus diagram, which consists of an abacus with one infinite horizontal runner; this runner has positions marked with the integers increasing from left to right, and has a bead at the position corresponding to each element of the beta-set.  For example, the following diagram depicts the beta-set $B^1(\la)$, where $\la$ is the partition $(4^2,2)$ in the example above.

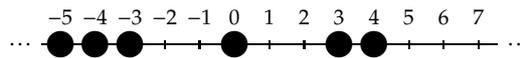
\[
\begin{smallmatrix}
&-5&-4&-3&-2&-1&0&1&2&3&4&5&6&7&\\
\hhd&\hbd&\hbd&\hbd&\hnb&\hnb&\hbd&\hnb&\hnb&\hbd&\hbd&\hnb&\hnb&\hnb&\hhd
\end{smallmatrix}
\]

\subsubs{Row removal and conjugate partitions}

For any partition $\la$, we denote by $\ol\la$ the partition $(\la_2,\la_3,\dots)$ obtained by removing the first part from $\la$.  We denote by $\la'$ the \emph{conjugate} partition to $\la$, given by
\[\la'_i=\left|\rset{j}{i\ls\la_j}\right|.\]

\subs{Specht modules}

For every partition $\la$ of $n$, there is an $\aaa n$-module $S^\la$ called a \emph{Specht module}.  In the case where $\aaa n$ is semisimple, the Specht modules are irreducible and afford all irreducible representations of $\aaa n$.  In general, the irreducible representations of $\aaa n$ can be obtained as quotients of certain Specht modules.  A similar situation applies for the algebra $\bbb n$, except that here the Specht modules are indexed by bipartitions of $n$.  The definition of Specht modules is given in \S\ref{murphy}.

Although not all Specht modules are irreducible, some are, and this paper is a contribution towards the classification of irreducible Specht modules.  For $\aaa n$, this question has been studied in a series of papers, and answered in all cases except where $q=\mo$.  In the present paper, we begin the study of irreducible Specht modules for $\bbb n$.

\subs{Simple modules and regular bipartitions}\label{conve}

We now briefly address the classification of simple modules for $\aaa n$ and $\bbb n$.  There are various conventions for (bi)partitions and Specht modules, which reflect established conventions in the representation theory of algebraic groups, quantum groups, and symmetric groups.  In practice it is very easy to translate between different conventions, but for clarity we set out here which convention we use.  This is the \emph{regular} convention: we say that a partition is \emph{$e$-regular} if either $e=\infty$ or $e<\infty$ and there is no $i$ such that $\la_i=\la_{i+e-1}>0$.  Then for every $e$-regular partition $\la$, the Specht module $S^\la$ for $\aaa n$ has an irreducible cosocle $D^\la$, and the modules $D^\la$ give all the irreducible representations of $\aaa n$.  Moreover, the composition multiplicity $[S^\la:D^\mu]$ equals $1$ for $\la=\mu$, and is zero unless $\mu\dom\la$.

For $\bbb n$, the simple modules are labelled by a certain class of bipartitions called \emph{regular} bipartitions in \cite{bk} (or \emph{conjugate Kleshchev} bipartitions in \cite{mfwt2}).  The Specht module $S^{(\la,\mu)}$ has an irreducible quotient $D^{(\la,\mu)}$ for each regular $(\la,\mu)$, and these modules afford all the irreducible representations of $\bbb n$.  The composition multiplicity $[S^{(\la,\mu)}:D^{(\nu,\xi)}]$ is zero unless $(\nu,\xi)\dom(\la,\mu)$, and equals $1$ if $(\nu,\xi)=(\la,\mu)$.

The definition of regular bipartitions is quite complicated in general, and depends on $q$ and $Q$; the original recursive definition \cite{arma} derives from the theory of crystals.  In the case $e=\infty$, it is easy to derive a non-recursive definition.

\begin{propn}\label{infreg}
Suppose $e=\infty$, and $(\la,\mu)$ is a bipartition.  Then $(\la,\mu)$ is regular if and only if one of the following holds:
\begin{itemize}
\item
$\mq$ is not a power of $q$;
\item
$\mq=q^r$ for some $r\gs0$, and $\la_i\gs \mu_i-r$ for all $i\gs1$;
\item
$\mq=q^r$ for some $r\ls0$, and $\la_i\gs\mu_{i+r}$ for all $i\gs1$.
\end{itemize}
\end{propn}

\begin{pf}
This is a simple exercise using the crystal-theoretic definition of a regular bipartition.  It also follows from Theorem \ref{decomp} below.
\end{pf}

A simple description of regular bipartitions is also easy to obtain in the case $e=2$ \cite[Proposition 4.11]{msimp}.  For the cases where $2<e<\infty$, a quick (though still recursive) characterisation of regular bipartitions has been found by Ariki, Kreiman and Tsuchioka \cite{akt}.  

We remark that the problem of classifying the irreducible Specht modules labelled by regular bipartitions (that is, Specht modules $S^{(\la,\mu)}=D^{(\la,\mu)}$) has been solved by James and Mathas, using their higher-level version of the Jantzen sum formula \cite[Theorem 4.7(iii)]{jmjs2}.

\subs{Conjugation and duality}\label{duality}

The anti-automorphism of $\bbb n$ which fixes each of the generators $T_0,\dots,T_{n-1}$ allows one to define a duality $M\mapsto M^\ast$ on $\bbb n$-modules in the natural way.  In order to understand the effect of this map on Specht modules, we need also to consider the automorphism $\theta$ of $\bbb n$ given by
\[T_0\longmapsto -QT_0^{-1},\qquad\qquad T_i\longmapsto -qT_i^{-1}\quad (i=1,\dots,n-1).\]
By twisting the action of $\bbb n$ by $\theta$, we obtain a self-equivalence $M\mapsto M^\theta$ on the category of $\bbb n$-modules.  The composition of this map with the duality $M\mapsto M^\ast$ will be denoted $M\mapsto M^\circledast$.  All these functors are self-inverse.

Now we have the following result; recall that $\la'$ denotes the conjugate of a partition $\la$.

\begin{thm}\label{duali}
Suppose $(\la,\mu)$ is a bipartition.
\begin{enumerate}
\item
$(S^{(\la,\mu)})^\circledast\cong S^{(\mu',\la')}$.
\item
If $(\la,\mu)$ is regular, then $(D^{(\la,\mu)})^\ast\cong D^{(\la,\mu)}$.
\end{enumerate}
\end{thm}

\begin{pfenum}
\item
This is essentially the result of \cite[Corollary 5.7]{tilt}; see \cite[Proposition 3.1]{mfwt2} for more details.
\item
This follows from the theory of cellular bases: the anti-automorphism $T_i\mapsto T_i$ is the automorphism implicit in the definition of cellular bases, under which simple modules are self-dual \cite[Exercise 2.7(iii)]{math}.
\end{pfenum}

As a consequence, we get the following result, which is relevant to the
present paper.

\begin{cory}\label{conj}
Suppose $(\la,\mu)$ is a bipartition.  Then $S^{(\la,\mu)}$ is reducible if and only if $S^{(\mu',\la')}$ is.
\end{cory}

Almost identical results hold in the simpler case of $\aaa n$; defining all the functors in the same way (ignoring the generator $T_0$), one gets $(S^\la)^\circledast\cong S^{\la'}$ for any partition $\la$, and $(D^\la)^\ast\cong D^\la$ if $\la$ is $e$-regular.

\subs{The Murphy basis}\label{murphy}

Here we give the definition of Specht modules for $\aaa n$ and $\bbb n$; this is based on the \emph{Murphy basis}, given in \cite{djm}, though we modify the definition in accordance with the conventions outlined in \S\ref{conve}.

Let $t_1,\dots,t_{n-1}$ be the Coxeter generators of the symmetric group $\sss n$.  Given $w\in\sss n$, let $t_{i_1}\dots t_{i_l}$ be any reduced expression for $w$, and define $T_w=T_{i_1}\dots T_{i_l}$; By Matsumoto's Theorem and the defining relations for $\aaa n$, this definition does not depend on the choice of reduced expression.

Define $\sss{\la'}$ to be the Young subgroup of $\sss n$ corresponding to $\la'$; that is, the naturally embedded subgroup $\sss{\la'_1}\times\sss{\la'_2}\times\cdots$.  Now define
\[x_\la=\sum_{w\in\sss{\la'}}T_w.\]
Let $\ol N^\la$ be the two-sided ideal of $\aaa n$ generated by $\lset{x_\mu}{\mu\domsby\la}$, and let $M^\la$ be the left ideal of $\aaa n$ generated by $x_\la$.  Then the Specht module is
\[
S^\la:=\frac{M^\la}{M^\la\cap \ol N^{\la}}.
\]

Now we consider bipartitions.  Given $0\ls a\ls n$, define
\[u_a=\prod_{i=1}^a\left(q^{i-1}+T_{i-1}T_{i-2}\dots T_1T_0T_1\dots
T_{i-1}\right).\]
Given a bipartition $(\la,\mu)$, define $\sss{(\mu',\la')}$ to be the Young subgroup $\sss{\mu'}\times\sss{\la'}$ of $\sss n$, and set
\[x_{(\la,\mu)}=u_{|\mu|}\sum_{w\in\sss{(\mu',\la')}}T_w.\]
As in type $A$, define $\ol N^{(\la,\mu)}$ to be the two-sided ideal of $\bbb n$ generated by $\rset{x_{(\nu,\xi)}}{(\nu,\xi)\domsby(\la,\mu)}$, and define $M^{(\la,\mu)}$ to be the left ideal generated by $x_{(\la,\mu)}$.  Now the
Specht module is
\[
S^{(\la,\mu)}:=\frac{M^{(\la,\mu)}}{M^\la\cap\ol N^{(\la,\mu)}}.
\]

We use this definition of the Specht modules to prove the following statement; this does not seem to appear in this form in the literature, but it admits a simple interpretation in terms of modules for the affine Hecke algebra of type $A$.

\begin{propn}\label{firstempty}
Suppose $\la$ is a partition of $n$.  Then the Specht module $S^{(\la,\varnothing)}$ (or $S^{(\varnothing,\la)}$) for $\bbb n$ is irreducible if and only if the Specht module $S^\la$ for $\aaa n$ is irreducible.
\end{propn}

\begin{pf}
We consider the Specht module $S^{(\la,\varnothing)}$; the result for $S^{(\varnothing,\la)}$ will then follow from the results in \S\ref{duality}.  We use the notation relating to the Murphy basis from the previous section.

Let $N$ be the two-sided ideal of $\bbb n$ generated by $u_1=1+T_0$.  Then $\bbb n/N$ is naturally isomorphic to $\aaa n$; in particular, $\bbb n/N\cong\aaa n$ as left $\aaa n$-modules.  Under this isomorphism, any $x_{(\xi,\mu)}$ with $|\mu|>0$ maps to zero, since it contains $u_1$ as a factor.  On the other hand, $x_{(\xi,\varnothing)}$ maps to $x_\xi$ for any $\xi$.  Hence $M_{(\la,\varnothing)}$ maps to $M_\la$, and $N_{\domsby(\la,\varnothing)}$ maps to $N_{\domsby\la}$.  As a consequence, we see that $S^{(\la,\varnothing)}\cong S^\la$ as $\aaa n$-modules.  Clearly, $T_0$ acts as the scalar $\mo$ on $\bbb n/N$, and so $S^{(\la,\varnothing)}$ and $S^\la$ have identical submodule structures.
\end{pf}

\subs{Morita equivalence}\label{morita}

We now cite a result of Dipper and James which will allow us to reduce the classification of irreducible Specht modules to the type $A$ question given certain assumptions about our parameters $q,Q$.

\begin{thm}\thmcite{dj}{Theorem 4.17}\label{djmorita}
Suppose $\mq$ is not a power of $q$.  Then $\bbb n$ is Morita equivalent to
\[\bigoplus_{m=0}^n\aaa m\otimes\aaa{n-m}.\]
\end{thm}

The Morita equivalence in Theorem \ref{djmorita} is constructed in such a way that the Specht module $S^{(\la,\mu)}$ corresponds to the product $S^\la\otimes S^\mu$ of Specht modules for Hecke algebras of type $A$. Hence we have the following result.

\begin{cory}\label{split}
Suppose $\mq$ is not a power of $q$ in $\bbf$, and $(\la,\mu)$ is a bipartition of $n$.  Then the Specht module $S^{(\la,\mu)}$ is irreducible if and only if the Specht modules $S^\la$ and $S^\mu$ in type $A$ are both irreducible.
\end{cory}

In fact, Theorem \ref{djmorita} has been generalised to higher levels by Dipper and Mathas \cite{dima}, so Corollary \ref{split} can be generalised to Ariki--Koike algebras.  In view of Corollary \ref{split}, we shall feel free to assume from now on that $Q=-q^r$ for $r\in\bbz$.

\subs{Decomposition maps}\label{decompmap}

In this section, we explain briefly some decomposition maps which will allow us to compare different Iwahori--Hecke algebras.  The framework for decomposition maps for Iwahori--Hecke algebras is explained in the article by Geck \cite{geck}.  Briefly, the idea is that one defines the Hecke algebra over an integral domain $A$ rather than a field.  By extending scalars to the field of fractions $K$ of $A$, one obtains the Hecke algebra $H_K$ over that field.  On the other hand, given a prime ideal $\fkp$ in $A$, one can extend scalars to the local ring $A_{\fkp}$, and then quotient by $\fkp$ to obtain the Hecke algebra $H_k$ algebra over the field $k=A/\fkp$.  By defining Specht modules over $A_{\fkp}$, one obtains a \emph{decomposition map}, which is a homomorphism between the Grothendieck groups of $H_K$ and $H_k$, sending the class of a Specht module for $H_K$ to the class of the corresponding Specht module for $H_k$.

Decomposition maps are very useful from our point of view, since (the class of) a simple module is mapped to a non-zero sum of simple modules (which essentially reflects how the simple module decomposes modulo $\fkp$).  Hence if a Specht module for $H_K$ is reducible, the corresponding Specht module for $H_k$ must be reducible.

There are two particular decomposition maps we shall use in this paper.  To explain these, we expand our notation, writing $\bbb n(\bbf,Q,q)$ for the Iwahori--Hecke algebra over $\bbf$ with parameters $Q,q$.

\begin{enumerate}
\item
Let $\hat q$ be an indeterminate over $\bbf$, and set $A=\bbf[\hat q]$.  Let $\fkp$ be the ideal of $A$ generated by $\hat q+1$; then there is a decomposition map between the Grothendieck groups of
\[\bbb n\big(\bbf(\hat q),-\hat q^r,\hat q\big)\qquad\text{and}\qquad\bbb n\big(\bbf,(\mo)^{r+1},\mo\big).\]
Note that the algebra on the left has $e=\infty$, while that on the right has $e=2$.
\item
Let $\hat Q$ be an indeterminate over $\bbf$, and set $A=\bbf[\hat Q]$.  For $r\in\{0,1\}$, let $\fkp$ be the ideal generated by $\hat Q+(\mo)^r$; then there is a decomposition map between the Grothendieck groups of
\[\bbb n\big(\bbf(\hat Q),\hat Q,\mo\big)\qquad\text{and}\qquad\bbb n\big(\bbf,(\mo)^{r+1},\mo\big).\]
Now both algebras have $e=2$, but in the algebra on the left we have $\hat Q\neq\pm1$, so we may apply the results in \S\ref{morita}.
\end{enumerate}

\subs{Blocks, induction and restriction}\label{indres}

For $n>1$, the algebra $\bbb{n-1}$ is naturally a subalgebra of $\bbb n$.  Here we briefly consider induction and restriction functors between the module categories of these algebras.

First we need to address the block structure of $\bbb n$.  Since the Specht modules are the cell modules for a cellular structure on $\bbb n$, it follows from \cite[Theorem 3.7]{grle} that the composition factors of any Specht module all lie in the same block of $\bbb n$.  So in order to describe the block structure, we just need to say when two Specht modules lie in the same block.  To do this, recall the definition of the residues of nodes of the Young diagram of a bipartition.

\begin{thm}\thmcite{lyma}{Theorem 2.11}
Suppose $(\la,\mu)$ and $(\nu,\xi)$ are bipartitions of $n$. Then $S^{(\la,\mu)}$ and $S^{(\nu,\xi)}$ lie in the same block of $\bbb n$ if and only if $(\la,\mu)$ and $(\nu,\xi)$ have the same number of nodes of residue $x$, for each $x\in\bbf$.
\end{thm}

As a result of this theorem, we may label a block of $\bbb n$ by a multiset of $n$ elements of $\bbf$.  Assuming for the rest of this subsection that $Q=-q^r$ for $r\in\bbz$, all these elements of $\bbf$ are powers of $q$.

Now we consider the natural functors $\ind$ and $\res$ of induction and restriction between $\bbb{n-k}$ and $\bbb n$.  Given a block $B$ of $\bbb{n-k}$ and a block $C$ of $\bbb n$, we have functors $\ind^C$ and $\res_B$ which act on modules lying in these blocks: $\ind^C$ is simply induction followed by projection onto $C$, and $\res_B$ is defined similarly.  In the case where the multiset corresponding to $C$ is obtained from the multiset corresponding to $B$ by adding $k$ copies of $q^i$ for some $i$, it is known that for any $M$ in $B$, $\ind^C(M)$ is a direct sum of $k!$ isomorphic modules; we write $f^{(k)}_i(M)$ for the isomorphism type of these modules.  We write $e^{(k)}_i(M)$ for the $B$-module obtained in the same way from $M$ in $C$.  Extending additively, $e^{(k)}_i$ and $f^{(k)}_i$ are defined on all $\bbb n$-modules, for all $n$.

The behaviour of these functors on Specht modules is quite well understood, via the following theorem.

\begin{thm}\thmlabel{The Branching Rule}
Suppose $(\la,\mu)$ is a bipartition.
\begin{itemize}
\item
\thmcite{mbranch}{Main Theorem} Let $(\la,\mu)^{+1},\dots,(\la,\mu)^{+s}$ be the bipartitions which can be obtained by adding $k$ nodes of residue $q^i$ to $(\la,\mu)$.  Then $f^{(k)}_i(S^{(\la,\mu)})$ has a filtration in which the factors are $S^{(\la,\mu)^{+1}},\dots,S^{(\la,\mu)^{+s}}$.
\item
\thmcite{arma}{Proposition 1.9} Let $(\la,\mu)^{-1},\dots,(\la,\mu)^{-t}$ be the bipartitions which can be obtained by removing $k$ nodes of residue $q^i$ from $(\la,\mu)$.  Then $e^{(k)}_i(S^{(\la,\mu)})$ has a filtration in which the factors are $S^{(\la,\mu)^{-1}},\dots,S^{(\la,\mu)^{-t}}$.
\end{itemize}
\end{thm}

In the context of irreducible Specht modules, this result sometimes allows us to reduce a bipartition to a smaller one.

\begin{propn}\label{branch}
Suppose $(\la,\mu)$ is a bipartition, $Q=-q^r$ for $r\in\bbz$ and $i\in\bbz$ is such that $(\la,\mu)$ has no addable $q^i$-nodes.  Let $(\la,\mu)^-_i$ be the bipartition obtained by removing all the removable $q^i$-nodes of $(\la,\mu)$.  Then $S^{(\la,\mu)}$ is irreducible if and only if $S^{(\la,\mu)^-_i}$ is.
\end{propn}

\begin{pf}
This is a straightforward consequence of the Branching Rule.  The reader who requires more detail should look at the proof of \cite[Lemma 3.3]{mfirred}.
\end{pf}

For example, suppose $q=Q=\mo$ and $(\la,\mu)=\left((4,3^2,1),(2,1)\right)$. 
The Young diagram $[\la,\mu]$, with the residues of all addables and removable
nodes marked, is as follows.
\Yboxdim{15pt}
\[\gyoung(;;;;\mo:\po::;;\mo:\po,;;;:\po:::;\mo:\po,;;;\po:::::\po,;\mo:\po,:\po)\]
We see that there are no addable nodes of residue $\mo$, and so $S^{(\la,\mu)}$ is irreducible if and only if $S^{(\la,\mu)^-_{\mo}}$ is, where $(\la,\mu)^-_{\mo}$ is the bipartition $\left((3^3),(1)\right)$ obtained by removing the removable nodes of residue $\mo$.

\vspace{\topsep}
We shall also need to consider the effect of the functor $e_i^{(k)}$ on simple modules, which is described by the \emph{modular branching rules}.  The first versions of these rules were proved for the symmetric groups by Kleshchev, and generalised to Hecke algebras of type $A$ by Brundan.  The generalisation to higher levels that we use here was proved by Ariki \cite[Theorem 6.1]{ari2}, building on the work of Grojnowski and Vazirani \cite{gv,g}.  We shall use only the following result, which combines a special case of the modular branching rules with the block classification.

\begin{propn}\label{modbranch}
Suppose $Q=-q^r$ for $r\in\bbz$.  Suppose $(\nu,\xi)\neq(\pi,\rho)$ are regular bipartitions of $n$.  Suppose $i\in\bbz$, and let $k$ be maximal such that $e^{(k)}_i(D^{(\nu,\xi)})\neq0$. Then:
\begin{enumerate}
\item\label{ressimp}
$e^{(k)}_i(D^{(\nu,\xi)})$ is simple, and the regular bipartition labelling this simple module does not depend on the characteristic of $\bbf$;
\item\label{otherreg}
if $e=\infty$, the composition multiplicity $\left[e^{(k)}_i(D^{(\pi,\rho)}):e^{(k)}_i(D^{(\nu,\xi)})\right]$ is zero.
\end{enumerate}
\end{propn}

\begin{pf}[Sketch proof]
Part (\ref{ressimp}) is a central part of the modular branching rules.

For part (2), we may assume $k>0$.  We may also assume $D^{(\nu,\xi)}$ and $D^{(\pi,\rho)}$ lie in the same block of $\bbb n$, since otherwise $e^{(k)}_i(D^{(\pi,\rho)})$ and $e^{(k)}_i(D^{(\nu,\xi)})$ lie in distinct blocks of $\bbb{n-k}$.  The block classification implies that for two bipartitions labelling simple modules in the same block, the number of removable $q^i$-nodes minus the number of addable $q^i$-nodes is the same; this proved in \cite[Proposition 3.2]{weight}, and is a consequence of the relationship between the block classification and the weight space decomposition for a highest-weight module for a certain Kac--Moody algebra.  The assumption $e=\infty$ means that the total number of addable and removable $q^i$-nodes of a bipartition is at most $2$ (since there can be at most one in each component).  Since $k$ is at most the number of removable $q^i$-nodes of $(\nu,\xi)$, we have $k\ls2$.

If $k=2$, then $(\nu,\xi)$ must have two removable and no addable $q^i$-nodes, and hence the same is true of $(\pi,\rho)$.  Hence by the modular branching rules $e^{(2)}_i(D^{(\pi,\rho)})$ is a simple module different from $e^{(2)}_i(D^{(\nu,\xi)})$.

If $k=1$, then the numbers of addable and removable $q^i$-nodes of $(\pi,\rho)$ can be any of $(0,0)$, $(1,1)$ or $(0,1)$.  Now the modular branching rules imply that $e_i(D^{(\pi,\rho)})$ is either zero or a simple module different from $e_i(D^{(\nu,\xi)})$.
\end{pf}

\begin{rmk}
We remark that it is vital in part (2) of the above proposition that $e=\infty$ and that we are working in type $B$ (rather than with an Ariki--Koike algebra of higher level).  If either of these assumptions is dropped, then it is possible to have $k\gs3$, and the argument in the above proof fails to work.
\end{rmk}

\sect{The case $e=\infty$}\label{infsec}

In this section, we give the classification of irreducible Specht modules for $\bbb n$ in the case where $q$ is not a root of unity.  This will be fairly straightforward, because the decomposition numbers can be given explicitly in this case.  In fact, this is not particularly new, but we give a detailed account below in order to make this paper reasonably self-contained and to express the main result in a combinatorial way which is suitable for our purposes.

When $\nchar(\bbf)=\infty$, Ariki's Theorem \cite{ari} implies that the decomposition numbers for $\bbb n$ are given by specialising the coefficients in canonical basis elements for a certain highest-weight $U_v(\mathfrak{sl}_\infty)$-module $V$; these coefficients are computed explicitly in \cite{lemi}.  To complete the description of the decomposition numbers, it suffices to show that they are unchanged when the assumption on characteristic is dropped.  We give a proof of this fact below using the results of \S\ref{indres}, though results in the literature yield a proof which essentially boils down to the same argument: as observed by Brundan \cite[Lemma 3.19]{br}, all the strings in the crystal of $V$ have length at most $2$, and this implies that each canonical basis element can be written in the form $f_{i_1}\dots f_{i_r}u_\Lambda$, where $u_\Lambda$ is a highest-weight vector in $V$ \cite[Proposition 4]{lemi}; as observed by Ariki and Mathas \cite[Corollary 3.7]{arma2}, this guarantees that the decomposition numbers are independent of characteristic.  This result is also obtained in \cite[\S9]{bs3}.

In order to give the formula for decomposition numbers, we have to set up some combinatorics.  This is all taken from \cite[\S2]{lemi}, but modified to suit our needs and our conventions (bearing in mind \S\ref{conve}).  Closely related combinatorics (in the form of cup and cap diagrams) is introduced in \cite[\S2]{bs1}.

\subs{$\pm$-sequences}

We define a \emph{$\pm$-sequence} to be a finite string $s_1\dots s_n$ of signs $\pm$.  When writing a $\pm$-sequence, we may group together consecutive equal terms with a superscript.  We say that a $\pm$-sequence $s_1\dots s_n$ is \emph{\domi} if either
\begin{itemize}
\item
every initial segment $s_1\dots s_m$ contains at least as many $-$ signs as $+$ signs, or
\item
every final segment $s_m\dots s_n$ contains at least as many $+$ signs as $-$ signs.
\end{itemize}

Note that the first condition can only hold if there are at least as many $-$ signs as $+$ signs in the whole sequence, and the second condition can only hold if there are at least as many $+$ signs as $-$ signs altogether.  If there are the same number of $+$ signs and $-$ signs (we say that the sequence is \emph{balanced} in this case), the two conditions are equivalent.

In this paper, we define an \emph{involution} of $\{1,\dots,n\}$ to be a permutation $\iota$ of $\{1,\dots,n\}$ such that $\iota^2$ is the identity.  If $s=s_1\dots s_n$ is a \domi{} $\pm$-sequence, then it is easy to see that there is a unique involution $\iota_s$ of $\{1,\dots,n\}$ with the following properties:
\begin{itemize}
\item
if $i$ and $j$ are fixed by $\iota_s$, then $s_i=s_j$;
\item
if $\iota_s(i)>i$, then $s_i=-$ and $s_{\iota(i)}=+$;
\item
there do not exist $i<j<k$ such that $\iota_s(i)=k$, $\iota_s(j)=j$;
\item
there do not exist $i<j<k<l$ such that $\iota_s(i)=k$, $\iota_s(j)=l$.
\end{itemize}
We define an involution $\iota$ of $\{1,\dots,n\}$ to be \emph{compatible} with $s$ if for each $i$ we have $\iota(i)\in\{i,\iota_s(i)\}$.  Given an involution $\iota$ of $\{1,\dots,n\}$, we write $s^\iota$ for the $\pm$-sequence $s_{\iota(1)}\dots s_{\iota(n)}$.

(Note that the involution $\iota_s$ corresponds roughly to the function $\psi$ in \cite[\S2.5]{lemi}, and the sequences $s^\iota$ for compatible $\iota$ correspond to elements of the set $\mathcal{C}(S)$.)

\subs{Decomposition numbers in the case $e=\infty$}

Now we can describe the decomposition numbers for $\bbb n$ when $e=\infty$.  Following \S\ref{morita}, we assume $Q=-q^r$ for $r\in\bbz$.  Suppose $(\nu,\xi)$ is a bipartition.  Choose any $i\in\bbz$, and construct the beta-sets $B_1=B^{r+i}(\nu)$, $B_2=B^{i}(\xi)$.  (These are infinite versions of the sets $\beta,\gamma$ used in \cite{lemi}.)  Let $B$ be the symmetric difference
\[
B=B_1\syd B_2:=\left(B_1\cup B_2\right)\setminus\left(B_1\cap B_2\right).
\]
Note that $B$ is finite, being the symmetric difference of two beta-sets.  We write the elements of $B$ as $b_1<\dots<b_n$, and define the \emph{\infs r} of $(\nu,\xi)$ to be the $\pm$-sequence $s=s_1\dots s_n$, where
\[
s_i=\begin{cases}
+&(b_i\in B_1)\\
-&(b_i\in B_2).
\end{cases}
\]
Note that this sequence is independent of the choice of $i$.

Given any involution $\iota$ of $\{1,\dots,n\}$, define
\begin{align*}
B_1^\iota &= \left(B_1\cap B_2\right)\cup\lset{b_{\iota(i)}}{b_i\in B_1\setminus B_2},\\
B_2^\iota &= \left(B_1\cap B_2\right)\cup\lset{b_{\iota(i)}}{b_i\in B_2\setminus B_1}.
\end{align*}
$B_1^\iota$ and $B_2^\iota$ are both beta-sets, and so define a pair of partitions; we denote the corresponding bipartition $(\nu,\xi)^\iota$.

Now we have the following result.

\begin{thm}\label{decomp}
Suppose $e=\infty$ and $(\nu,\xi)$ is a bipartition, and define the set $B$ and the sequence $s=s_1\dots s_n$ as above.  Then $(\nu,\xi)$ is a regular bipartition if and only if $s$ is dominant.  If this is the case, then $D^{(\nu,\xi)}$ occurs once as a composition factor of $S^{(\nu,\xi)^\iota}$ for every involution $\iota$ of $\{1,\dots,n\}$ compatible with $s$, and does not occur as a composition factor of any other Specht module.
\end{thm}

\begin{pf}
By Ariki's Theorem \cite[Theorem 4.4]{ari}, the decomposition numbers for $\bbb n$ in the case where $\bbf$ has infinite characteristic can be obtained by specialising at $v=1$ the canonical basis coefficients for the irreducible highest-weight module $V(\La_r+\La_0)$ for the quantum algebra $U_v(\mathfrak{sl}_\infty)$.  These canonical basis coefficients are given explicitly in \cite[Theorem 3]{lemi}, so the result follows from a simple translation of notation.

To complete the proof, it remains to show that the decomposition numbers for $\bbb n$ are independent of the underlying characteristic when $e=\infty$.  
We prove this by induction on $n$.  Consider the decomposition number $[S^{(\la,\mu)}:D^{(\nu,\xi)}]$, where $(\la,\mu)$ and $(\nu,\xi)$ are bipartitions of $n$ with $(\nu,\xi)$ $e$-regular.  This decomposition number is zero if $S^{(\la,\mu)}$ and $D^{(\nu,\xi)}$ are not in the same block of $\bbb n$ (and the condition for them to lie in the same block is independent of the characteristic), so we assume they do lie in the same block.

The restriction of $D^{(\nu,\xi)}$ to $\bbb{n-1}$ is non-zero, so there must be some $i\in\bbz$ such that $e_i(D^{(\nu,\xi)})\neq0$.  Choose such an $i$, and let $k$ be maximal such that $e^{(k)}_i(D^{(\nu,\xi)})\neq0$.  By the Branching Rule and by induction, the composition factors of $e_i^{(k)}(S^{(\la,\mu)})$ are independent of the underlying characteristic.  If $D^{(\pi,\rho)}$ is a composition factor of $S^{(\la,\mu)}$ other than $D^{(\nu,\xi)}$, then by Proposition \ref{modbranch}(\ref{otherreg}) the composition multiplicity $\big[e^{(k)}_i(D^{(\pi,\rho)}):e^{(k)}_i(D^{(\nu,\xi)})\big]$ is zero.  Hence, since $e_i^{(k)}$ is an exact functor, the decomposition number $[S^{(\la,\mu)}:D^{(\nu,\xi)}]$ equals the multiplicity $\big[e^{(k)}_i(S^{(\la,\mu)}):e^{(k)}_i(D^{(\nu,\xi)})\big]$, and so by induction and Proposition \ref{modbranch}(\ref{ressimp}) is independent of the characteristic.
\end{pf}

\begin{eg}
Suppose $r=1$ and $(\nu,\xi)=\big((4^2,3^3),(4^2,1)\big)$.  Then
\begin{align*}
B^1(\nu)&=\{\dots,-7,-6,-5,-1,0,1,3,4\},\\
B^0(\xi)&=\{\dots,-7,-6,-5,-4,-2,2,3\},
\end{align*}
so that $B=\{-4,-2,-1,0,1,2,4\}$ and $s={-}{-}{+}{+}{+}{-}{+}$.  This is easily seen from the following diagram, in which we give the abacus diagram for $B^1(\nu)$ above that for $B^0(\xi)$; each bead directly above an empty position contributes a $+$ to $s$, and each bead below an empty position contributes a $-$.
\medskip
\setlength\unitlength{0.9pt}
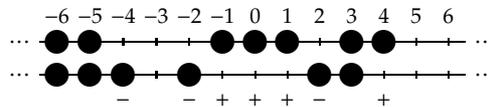
\[
\begin{smallmatrix}
&-6&-5&-4&-3&-2&-1&0&1&2&3&4&5&6&\\
\hhd&\hbd&\hbd&\hnb&\hnb&\hnb&\hbd&\hbd&\hbd&\hnb&\hbd&\hbd&\hnb&\hnb&\hhd\\
\hhd&\hbd&\hbd&\hbd&\hnb&\hbd&\hnb&\hnb&\hnb&\hbd&\hbd&\hnb&\hnb&\hnb&\hhd\\
&&&-&&-&+&+&+&-&&+&&&&
\end{smallmatrix}
\]
\medskip
$s$ is a dominant sequence, so $(\nu,\xi)$ is a regular bipartition.  The involution $\iota_s$ is $(1,4)(2,3)(6,7)$.  Hence there are eight Specht modules containing the simple module $D^{(\nu,\xi)}$; the labelling partitions with the corresponding abacus diagrams are given below.
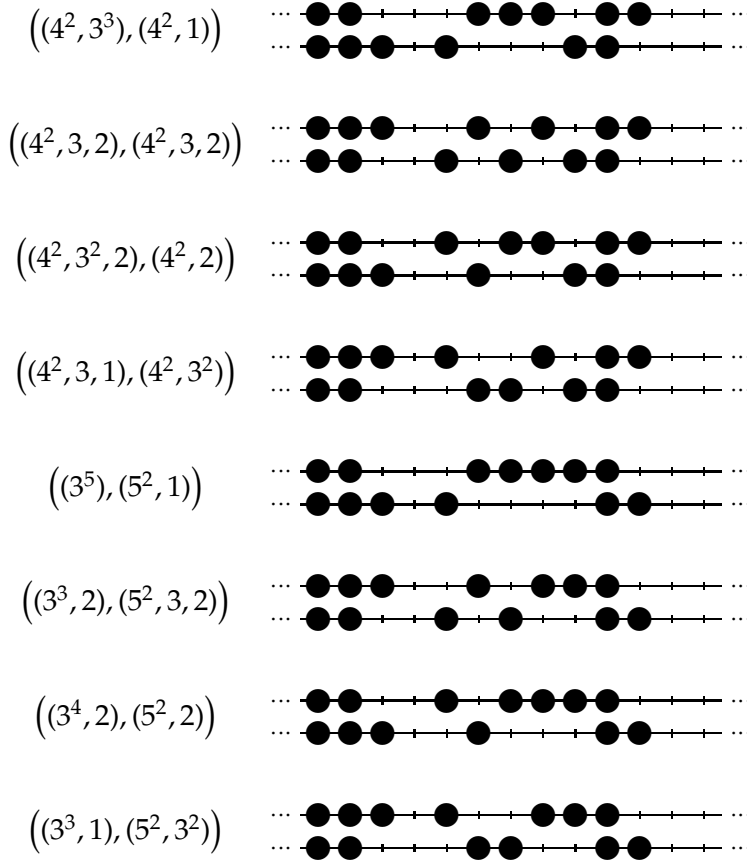
\newlength\coll
\setlength\coll{100pt}
\newlength\skp
\setlength\skp{16pt}
{\allowdisplaybreaks\begin{alignat*}2
&\hbox to \coll{\hfil $\big((4^2,3^3),(4^2,1)\big)$\hfil }&&
\begin{smallmatrix}
\hhd&\hbd&\hbd&\hnb&\hnb&\hnb&\hbd&\hbd&\hbd&\hnb&\hbd&\hbd&\hnb&\hnb&\hhd\\
\hhd&\hbd&\hbd&\hbd&\hnb&\hbd&\hnb&\hnb&\hnb&\hbd&\hbd&\hnb&\hnb&\hnb&\hhd\\
\end{smallmatrix}
\\[\skp]
&\hbox to \coll{\hfil $\big((4^2,3,2),(4^2,3,2)\big)$\hfil }&&
\begin{smallmatrix}
\hhd&\hbd&\hbd&\hbd&\hnb&\hnb&\hbd&\hnb&\hbd&\hnb&\hbd&\hbd&\hnb&\hnb&\hhd\\
\hhd&\hbd&\hbd&\hnb&\hnb&\hbd&\hnb&\hbd&\hnb&\hbd&\hbd&\hnb&\hnb&\hnb&\hhd\\
\end{smallmatrix}
\\[\skp]
&\hbox to \coll{\hfil $\big((4^2,3^2,2),(4^2,2)\big)$\hfil }&&
\begin{smallmatrix}
\hhd&\hbd&\hbd&\hnb&\hnb&\hbd&\hnb&\hbd&\hbd&\hnb&\hbd&\hbd&\hnb&\hnb&\hhd\\
\hhd&\hbd&\hbd&\hbd&\hnb&\hnb&\hbd&\hnb&\hnb&\hbd&\hbd&\hnb&\hnb&\hnb&\hhd\\
\end{smallmatrix}
\\[\skp]
&\hbox to \coll{\hfil $\big((4^2,3,1),(4^2,3^2)\big)$\hfil }&&
\begin{smallmatrix}
\hhd&\hbd&\hbd&\hbd&\hnb&\hbd&\hnb&\hnb&\hbd&\hnb&\hbd&\hbd&\hnb&\hnb&\hhd\\
\hhd&\hbd&\hbd&\hnb&\hnb&\hnb&\hbd&\hbd&\hnb&\hbd&\hbd&\hnb&\hnb&\hnb&\hhd\\
\end{smallmatrix}
\\[\skp]
&\hbox to \coll{\hfil $\big((3^5),(5^2,1)\big)$\hfil }&&
\begin{smallmatrix}
\hhd&\hbd&\hbd&\hnb&\hnb&\hnb&\hbd&\hbd&\hbd&\hbd&\hbd&\hnb&\hnb&\hnb&\hhd\\
\hhd&\hbd&\hbd&\hbd&\hnb&\hbd&\hnb&\hnb&\hnb&\hnb&\hbd&\hbd&\hnb&\hnb&\hhd\\
\end{smallmatrix}
\\[\skp]
&\hbox to \coll{\hfil $\big((3^3,2),(5^2,3,2)\big)$\hfil }&&
\begin{smallmatrix}
\hhd&\hbd&\hbd&\hbd&\hnb&\hnb&\hbd&\hnb&\hbd&\hbd&\hbd&\hnb&\hnb&\hnb&\hhd\\
\hhd&\hbd&\hbd&\hnb&\hnb&\hbd&\hnb&\hbd&\hnb&\hnb&\hbd&\hbd&\hnb&\hnb&\hhd\\
\end{smallmatrix}
\\[\skp]
&\hbox to \coll{\hfil $\big((3^4,2),(5^2,2)\big)$\hfil }&&
\begin{smallmatrix}
\hhd&\hbd&\hbd&\hnb&\hnb&\hbd&\hnb&\hbd&\hbd&\hbd&\hbd&\hnb&\hnb&\hnb&\hhd\\
\hhd&\hbd&\hbd&\hbd&\hnb&\hnb&\hbd&\hnb&\hnb&\hnb&\hbd&\hbd&\hnb&\hnb&\hhd\\
\end{smallmatrix}
\\[\skp]
&\hbox to \coll{\hfil $\big((3^3,1),(5^2,3^2)\big)$\hfil }&&
\begin{smallmatrix}
\hhd&\hbd&\hbd&\hbd&\hnb&\hbd&\hnb&\hnb&\hbd&\hbd&\hbd&\hnb&\hnb&\hnb&\hhd\\
\hhd&\hbd&\hbd&\hnb&\hnb&\hnb&\hbd&\hbd&\hnb&\hnb&\hbd&\hbd&\hnb&\hnb&\hhd\\
\end{smallmatrix}
\\
\end{alignat*}}
\end{eg}

We shall use Theorem \ref{decomp} to classify the irreducible Specht modules.  Given a bipartition $(\la,\mu)$, the Specht module $S^{(\la,\mu)}$ is irreducible if and only if there is a unique regular bipartition $(\nu,\xi)$ such that $(\la,\mu)$ occurs as $(\nu,\xi)^\iota$ for a compatible involution $\iota$.  It turns out that this condition is quite easy to express in terms of the \infs r{} of $(\la,\mu)$.

Given a dominant $\pm$-sequence $s$ and a compatible involution $\iota$ for $s$, we will say that the pair $(s,\iota)$ is \emph{suitable} for the sequence $s^\iota$.

\begin{propn}\label{onlyone}
Suppose $t$ is a $\pm$-sequence of the form ${+}^a{-}^b{+}^c$ or ${-}^a{+}^b{-}^c$, with $a+c\ls b$.  Then there is exactly one pair $(s,\iota)$ suitable for $t$.
\end{propn}

\begin{pf}
We assume $t$ has the form ${+}^a{-}^b{+}^c$; the proof in the other case is similar.  The pair $({-}^a{+}^a{-}^{b-a}{+}^c,\kappa)$, where
\[\kappa:i\longmapsto
\begin{cases}
2a+1-i&(i\ls 2a)\\
i&(i>2a)
\end{cases}
\]
is certainly suitable for $t$.  Now suppose $(s,\iota)$ is suitable for $t$.  Then we must have $s_{a+b+1}=\dots=s_{a+b+c}=+$, and $\iota$ must fix all the points $a+b+1,\dots,a+b+c$.  Since $t$ contains at least as many $-$ signs as $+$ signs, $s$ does too, and so since $s$ is dominant we must have $s_1=-$.  In fact, we claim that $s_i=-$ for all $1\ls i\ls a$.  If not, let $i$ be minimal such that $s_i=+$.  Then since $s_{i-1}=-$, we must have $\iota_s(i)=i-1$, and so $t_{i-1},t_i$ equal $+,-$ in some order, which is not the case.

So $s_1=\dots=s_a=-$, and this means that $1,\dots,a$ are all moved by $\iota$.  If we write $j_i=\iota(i)$ for $i=1,\dots,a$, then from the definition of $\iota_s$ we have $a+1\ls j_a<j_{a-1}<\dots<j_1\ls a+b$.  Now we find that $s_{j_a}=+$ and $s_{j_a-1}=-$, and this implies that $\iota_s(j_a)=j_a-1$.  But we have $\iota_s(j_a)=\iota(j_a)=a$, so that $j_a=a+1$.  And we claim that $j_i=2a+1-i$ for all $i$: if not, let $i$ be maximal such that $j_i>2a+1-i$; then we have $s_{j_i}=+$ and $s_{j_i-1}=-$, so that $\iota_s(j_i)=j_i-1$.  But this contradicts the fact that $\iota(i_j)=j$.

And now we have $s={-}^a{+}^a{-}^{b-a}{+}^c$, and $\iota=\kappa$.
\end{pf}

We want to prove a converse to Proposition \ref{onlyone}.  For the inductive step, we shall need the following lemma.

\begin{lemma}\label{plusminus}
Suppose $s=s_1\dots s_n$ is a dominant $\pm$-sequence, and $1\ls i<n$ is such that $s_i=-$, $s_{i+1}=+$.  Define $\hat s=s_1\dots s_{i-1}s_{i+2}\dots s_n$, and for any $j\in\{1,\dots,n-2\}$ write
\[
j''=
\begin{cases}
j&(j<i)\\
j+2&(j\gs i).
\end{cases}
\]
Then:
\begin{enumerate}
\item
$\hat s$ is a dominant sequence;
\item
if $\hat\iota$ is a compatible involution for $\hat s$, then the involutions $\iota^1,\iota^2$ of $\{1,\dots,n\}$ defined by
\[
\iota^1(j)=
\begin{cases}
{\hat\iota(j)}''&(j<i)\\
i&(j=i)\\
i+1&(j=i+1)\\
{\hat\iota(j-2)}''&(j>i+1)
\end{cases}\qquad
\iota^2(j)=
\begin{cases}
{\hat\iota(j)}''&(j<i)\\
i+1&(j=i)\\
i&(j=i+1)\\
{\hat\iota(j-2)}''&(j>i+1)
\end{cases}
\]
are compatible with $s$.
\end{enumerate}
\end{lemma}

\begin{pf}
It is immediate from the definition of a dominant sequence that $\hat s$ is dominant.  Moreover, it is easy to compare the involutions $\iota_s$ and $\iota_{\hat s}$: $\iota_s$ is given by
\[
j\longmapsto
\begin{cases}
{\iota_{\hat s}(j)}''&(j<i)\\
i+1&(j=i)\\
i&(j=i+1)\\
{\iota_{\hat s}(j-2)}''&(j>i+1).
\end{cases}
\]
Now the lemma follows.
\end{pf}

\begin{cory}\label{indcomp}
Suppose $t=t_1\dots t_n$ is a $\pm$-sequence with $t_i\neq t_{i+1}$, and $(u,\kappa)$ is a suitable pair for $\hat t=t_1\dots t_{i-1}t_{i+2}\dots t_n$.  Then there is a suitable pair $(s,\iota)$ for $t$ such that $s_j=u_j$ for $j\ls i-1$ and $s_j=u_{j-2}$ for $j\gs i+2$.
\end{cory}

\begin{pf}
Define $s$ by
\[
s_j=
\begin{cases}
u_j&(j\ls i-1)\\
-&(j=i)\\
+&(j=i+1)\\
u_{j-2}&(j\gs i+2).
\end{cases}
\]
Then apply Lemma \ref{plusminus} to $s$, taking $\iota=\iota^1$ if $t_i=-$, and $\iota=\iota^2$ if $t_i=+$.
\end{pf}

\begin{propn}\label{atleasttwo}
Suppose $t$ is a $\pm$-sequence which is not of the form ${-}^a{+}^b{-}^c$ or ${+}^a{-}^b{+}^c$ with $a,b,c$ non-negative integers and $a+c\ls b$.  Then there are at least two pairs $(s,\iota)$ suitable for $t$.
\end{propn}

\begin{pf}
Take any $i$ such that $t_i\neq t_{i+1}$, and let $\hat t=t_1\dots t_{i-1}t_{i+2}\dots t_n$.  If $\hat t$ does not have the form ${-}^a{+}^b{-}^c$ or ${+}^a{-}^b{+}^c$ with $a+c\ls b$ then by induction on $n$ the sequence $\hat t$ has at least two suitable pairs.  Applying Corollary \ref{indcomp}, we find that there are at least two pairs suitable for $t$.  So we may assume that \emph{for every} $i$ such that $t_i\neq t_{i+1}$ the sequence $\hat t$ has the form ${-}^a{+}^b{-}^c$ or ${+}^a{-}^b{+}^c$ with $a+c\ls b$.  This leaves only a few possibilities for $t$.
\begin{itemize}
\item
Suppose $t$ has the form ${+}^a{-}{+}^c$ with $a,c>0$.  Then the pair $(t,\id)$ is suitable for $t$, as is the pair $\left({+}^{a-1}{-}{+}^{c+1},(a,a+1)\right)$.  A similar argument applies in the case where $t={-}^a{+}{-}^c$ with $a,c>0$.
\item
Suppose $t$ has the form ${+}^a{-}{+}{-}^d$, where $a,d>0$.  If $a\ls d$, then we have two suitable pairs $\left({-}^{a+1}{+}^{a+1}{-}^{d-a},\iota\right)$ and $\left({-}^a{+}{-}{+}^a{-}^{d-a},\kappa\right)$, where
\begin{align*}
\iota:j&\longmapsto
\begin{cases}
2a+3-j&(j\ls a\text{ or }a+3\ls j\ls 2a+2)\\
j&(\text{otherwise}),
\end{cases}\\
\kappa:j&\longmapsto
\begin{cases}
2a+3-j&(j< a\text{ or }a+3< j\ls 2a+2)\\
j+1&(j=a\text{ or }a+2)\\
j-1&(j=a+1\text{ or }a+3)\\
j&(\text{otherwise}).
\end{cases}
\end{align*}
A similar argument applies when $a\gs d$.
\item
Suppose $t$ has the form ${-}^a{+}{-}{+}^d$ with $a,d>0$.  Then $(t,\id)$ is a suitable pair, as is $\left({-}^{a+1}{+}^{d+1},(a+b,a+b+1)\right)$.\qedhere
\end{itemize}
\end{pf}

Combining Propositions \ref{onlyone} and \ref{atleasttwo} with Theorem \ref{decomp} yields the main result of this section.

\begin{cory}\label{main0}
Suppose $Q=-q^r$ and $(\la,\mu)$ is a bipartition, and let $t$ be the \infs r{} of $(\la,\mu)$.  Then the Specht module $S^{(\la,\mu)}$ for $\bbb n$ is irreducible if and only if $t$ has the form ${-}^a{+}^b{-}^c$ or ${+}^a{-}^b{+}^c$ with $a,b,c$ non-negative integers and $a+c\ls b$.
\end{cory}

\begin{rmk}
In the case where $\nchar(\bbf)=\infty$, this result appears in \cite[Theorem 8.25]{bk1} in the context of representations of shifted Yangians; the combinatorial criterion on the signature first appears in \cite[Lemma 3.8]{lz}.
\end{rmk}

We shall define some terminology for use in Section \ref{e2sec}; let us say that a bipartition $(\la,\mu)$ is \emph{\infi r} if it satisfies the conditions of Corollary \ref{main0}, and \emph{\infr r} otherwise.

\sect{The case $e=2$}\label{e2sec}

\subs{The main result}

In this section, we consider the case where $e=2$; that is, $q=\mo$.  In this case we still do not have a classification of irreducible Specht modules for Iwahori--Hecke algebras of type $A$.  But our main result here is that the type $B$ problem can be reduced to the type $A$ problem.

To be precise, recall the definition of residues of nodes of Young diagrams.  Following the discussion in \S\ref{morita}, we assume that $Q=(\mo)^{r+1}$ for some integer $r$.  In fact, we assume that $r\in\{0,1\}$.  Recall that Proposition \ref{branch} sometimes allows us to reduce a bipartition to a smaller one.  If $i=\pm1$, let us say that a bipartition $(\la,\mu)$ is \emph{$i$-\dable} if $(\la,\mu)$ has at least one removable $i$-node and no addable $i$-nodes; say that $(\la,\mu)$ is \dable{} if it is either $1$-\dable{} or $\mo$-\dable.  Our main result in this section is the following.

\begin{thm}\label{main2}
Suppose $q=\mo$, $Q=\pm1$ and $(\la,\mu)$ is a bipartition which is not \dable{} and which has $\la\neq\varnothing\neq\mu$.  Then $S^{(\la,\mu)}$ is reducible.
\end{thm}

Combining this with Proposition \ref{branch} and Proposition \ref{firstempty}, we obtain the following corollary.  Let us say that a partition $\nu$ is $2$-irreducible if the Specht module $S^\nu$ for $\aaa n$ is irreducible.  (Note that whether $\nu$ is $2$-irreducible depends on the characteristic of $\bbf$; but the field $\bbf$ is fixed throughout this section.)

\begin{cory}\label{maincor}
Suppose $q=\mo$ and $Q=(\mo)^{r+1}$, and $(\la,\mu)$ is a bipartition.  Then $S^{(\la,\mu)}$ is irreducible if and only if there is $i=\pm1$ and a sequence $(\la,\mu)=(\la_0,\mu_0),\dots,(\la_s,\mu_s)$ of bipartitions such that:
\begin{itemize}
\item
for even $j$, $(\la_j,\mu_j)$ is $i$-\dable{} with $(\la_j,\mu_j)^-_i=(\la_{j+1},\mu_{j+1})$;
\item
for odd $j$, $(\la_j,\mu_j)$ is $(-i)$-\dable{} with $(\la_j,\mu_j)^-_{-i}=(\la_{j+1},\mu_{j+1})$;
\item
$(\la_s,\mu_s)$ has the form $(\nu,\varnothing)$ or $(\varnothing,\nu)$ for $\nu$ a $2$-irreducible partition.
\end{itemize}
\end{cory}

\begin{eg}
Suppose as in the example in \S\ref{indres} that $Q=\mo$ and $(\la,\mu)=\left((4,3^2,1),(2,1)\right)$.  Then in Corollary \ref{maincor} we take $i=\mo$ and
\begin{align*}
(\la_0,\mu_0)&=\left((4,3^2,1),(2,1)\right),\\
(\la_1,\mu_1)&=\left((3^3),(1)\right),\\
(\la_2,\mu_2)&=\left((3^2,2),\varnothing\right).
\end{align*}
The partition $(3^2,2)$ is known to be $2$-irreducible (since the characteristic of $\bbf$ is not $2$) and so $S^{(\la,\mu)}$ is irreducible.
\end{eg}

\subs{Decomposition maps}

Our main tool in proving Theorem \ref{main2} will be to use decomposition maps and our results for the case $e=\infty$.  Following the discussion in \S\ref{decompmap} and using the first decomposition map given in that section, we find that when $q=\mo$ and $Q=(\mo)^{r+1}$, the Specht module $S^{(\la,\mu)}$ for $\bbb n$ is irreducible only if $(\la,\mu)$ is \infi t \emph{for every $t$ of the same parity as $r$}.

With this in mind, we examine the condition for a bipartition to be \infi t in more detail.  Recall that if $\la$ is a partition then $\ol\la$ denotes the partition obtained by removing the first part from $\la$.

\begin{lemma}
Suppose $t\in\bbz$ and $(\la,\mu)$ is a bipartition, and let $s$ be the \infs t{} of $(\la,\mu)$.  Then the \infs t{} of $(\ol\la,\ol\mu)$ either equals $s$ or is obtained from $s$ by deleting the last $+$ and the last $-$.
\end{lemma}

\begin{pf}
Recall that $s$ may be constructed from the symmetric difference $\be t\la\syd\be0\mu$.  On the other hand, the \infs t{} of $(\ol\la,\ol\mu)$ may be constructed from $\be{t-1}{\ol\la}\syd\be{-1}{\ol\mu}$.  The set $\be{t-1}{\ol\la}$ is obtained from $\be t\la$ by removing the largest element $b$, while $\be{-1}{\ol\mu}$ is obtained from $\be0\mu$ by removing the largest element $c$.

If $b=c$, then clearly $\be t\la\syd\be0\mu=\be{t-1}{\ol\la}\syd\be{-1}{\ol\mu}$, so the \infs ts are the same.  If $c\neq b\in\be0\mu$, then the \infs t{} of $(\ol\la,\ol\mu)$ equals $s$: to obtain the \infs t of $(\ol\la,\ol\mu)$ from $s$, the last symbol (the $-$ corresponding to $c$) is removed, and a $-$ corresponding to $b$ is inserted; but since $b=\max(\be t\la)$, this insertion happens after the last $+$ in $s$, so the end result is that the \infs t is unchanged.  Similarly, if $b\neq c\in\be t\la$, then $s$ is the \infs t{} of $(\ol\la,\ol\mu)$.

The remaining possibility is that $b\notin\be0\mu$ and $c\notin\be t\la$.  In this case, the \infs t{} of $(\ol\la,\ol\mu)$ is obtained from $s$ by removing the last $+$ (corresponding to $b$) and the last $-$ (corresponding to $c$).
\end{pf}

As a consequence, we get the following.

\begin{lemma}\label{rowrem}
Suppose $t\in\bbz$ and  $(\la,\mu)$ is a bipartition which is \infi t.  Then $(\ol\la,\ol\mu)$ is \infi t.
\end{lemma}

\begin{pf}
Recall that $(\la,\mu)$ is \infi t{} if and only if the \infs t{} of $(\la,\mu)$ has the form ${+}^a{-}^b{+}^c$ or ${-}^a{+}^b{-}^c$ and with $a,c\gs0$ and $b\gs a+c$.  If this is the case, then the same is true of the same sequence with the last $+$ and the last $-$ removed.
\end{pf}

Now we come to the main part of the proof of Theorem \ref{main2}.

\begin{propn}\label{allirodd}
Suppose $q=\mo$, $Q=1$ and $(\la,\mu)$ is a bipartition which is \infi t{} for all odd $t$.  Then either
\begin{itemize}
\item
$\la$ or $\mu$ equals $\varnothing$, or
\item
$(\la,\mu)$ is \dable.
\end{itemize}
\end{propn}

\begin{pf}
First we note that $(\la,\mu)$ is \infi t{} if and only if $(\mu,\la)$ is \infi{-t}; this follows immediately from the definitions.  Also, $(\la,\mu)$ is \infi t{} if and only if $(\mu',\la')$ is (this can be seen by comparing signatures, or directly from Corollary \ref{conj}).  Hence, since $(\la,\mu)$ is \infi t{} for every odd value of $t$, the same is true of $(\mu,\la)$, $(\mu',\la')$ and $(\la',\mu')$; so we can interchange $\la$ and $\mu$, or replace them both with their conjugates, when it is convenient.

We shall suppose neither of the two conclusions is true, and derive a contradiction.  First suppose $\mu=(1)$.  Since $(\la,\mu)$ is not \dable, $\la$ has an addable node $(a,b)$ with $b-a$ odd.  But now we claim that $(\la,\mu)$ is \infr{a-b}.  We have $\be0\mu=\{\dots,-4,-3,-2,0\}$, while $\be{a-b}\la$ contains $-1$ but not $0$.  Moreover, since $\la\neq\varnothing$, there is either a positive integer in $\be{a-b}\la$, or an integer less than $-1$ not in $\be{a-b}\la$.  Hence the \infs{a-b} has the form ${-}^x{+}{-}{+}^y$ with at least one of $x$ and $y$ strictly positive.  So $(\la,\mu)$ is not \infi{a-b}.

So we can assume $|\mu|>1$.  Symmetrically, we can assume $|\la|>1$.  Now we claim that $(\la,\mu)$ is \infr{\la'_1-\mu_1-1}.  The set $\be{1-\mu_1}\mu$ has $0$ as its largest element, and (since $\mu\neq\varnothing$) does not contain every negative integer.  On the other hand, the set $\be{\la'_1}\la$ contains every negative integer, does not contain $0$ and (since $\la\neq\varnothing$) contains some positive integer.  So the \infs{(\la'_1-\mu_1-1)} has the form ${+}^x{-}{+}^y$, for some $x,y>0$, and $(\la,\mu)$ is \infr{\la'_1-\mu_1-1}.

Since $(\la,\mu)$ is \infr t for every odd $t$, we deduce that $\la'_1-\mu_1$ must be odd.  Symmetrically, $\mu'_1-\la_1$ must be odd.  This implies that $(\ol\la,\ol\mu)$ is not \dable.  To see this, suppose $(\ol\la,\ol\mu)$ has a removable node of residue $i=\pm1$.  Then $(\la,\mu)$ has a removable node of residue $-i$.  Since $(\la,\mu)$ is not \dable, this means that $(\la,\mu)$ has an addable node of residue $-i$.  Furthermore, this addable node can be chosen not to be $(1,\la_1+1)_1$ or $(1,\mu_1+1)_2$, because $(1,\la_1+1)_1$ has the same residue as the addable node $(\mu'_1+1,1)_2$, while $(1,\mu_1+1)_2$ has the same residue as $(\la'_1+1,1)_1$.  Hence $(\ol\la,\ol\mu)$ has an addable node of residue $i$, as required.

By induction, this implies that one of $\ol\la,\ol\mu$ equals $\varnothing$, i.e.\ $\la$ or $\mu$ has only one non-empty row.  Symmetrically, either $\la$ or $\mu$ has only one non-empty column.  Since neither $\la$ nor $\mu$ equals $(1)$, this means that (without loss of generality) $(\la,\mu)=\left((1^m),(n)\right)$ for some positive integers $m,n$.  The fact that $\mu'_1-\la_1$ and $\la'_1-\mu_1$ are odd means that $m$ and $n$ are both even.  But now one can check that the \infs{(m+n-1)} of $(\la,\mu)$ is ${+}^n{-}{+}^m$, so $(\la,\mu)$ is \infr{m+n-1}.
\end{pf}

The corresponding result for even values of $t$ is more complicated.

\begin{propn}\label{allireven}
Suppose $(\la,\mu)$ is a bipartition which is \infi t{} for all even $t$.  Then one of the following is true:
\begin{itemize}
\item
$\la$ or $\mu$ equals $\varnothing$;
\item
$(\la,\mu)$ is \dable;
\item
$(\la,\mu)=\left((2m),(1^{2n})\right)$ or $\left((1^{2n}),(2m)\right)$ for some $m,n$.
\end{itemize}
\end{propn}

\begin{pf}
We follow the proof of Proposition \ref{allirodd}.  As in that proof, we may interchange $\la$ and $\mu$ or replace them with their conjugates if necessary.  We begin by eliminating the case where $\la=(1)$; this follows exactly as before, replacing `odd' with `even'.

So we can assume $|\la|,|\mu|>1$.  We can also copy the proof of Proposition \ref{allirodd} to show that we must have $\mu'_1-\la_1$ and $\la'_1-\mu_1$ both even.  And as before this shows that $(\ol\la,\ol\mu)$ is not \dable.  Now by induction we find that either one of $\ol\la,\ol\mu$ equals $\varnothing$ or $(\ol\la,\ol\mu)$ has the form $\left((2m),(1^{2n})\right)$ or $\left((1^{2n}),(2m)\right)$.  A similar statement applies when we remove the first columns of $\la$ and $\mu$.  Together with what we have proved so far, the only remaining possibility is that $(\la,\mu)=\left((k),(1^l)\right)$ or $\left((1^l),(k)\right)$, where $k$ and $l$ have the same parity.  If $k$ and $l$ are both odd, then $(\la,\mu)$ is \dable, so $k$ and $l$ are even.
\end{pf}

\subs{Proof of Theorem \ref{main2}}

To complete the proof of Theorem \ref{main2}, we shall need the following result.

\begin{propn}\label{lacmu}
Suppose $e=2$, and $(\la,\mu)$ is a bipartition.  Then the Specht modules $S^{(\la,\mu)}$ and $S^{(\la',\mu)}$ have exactly the same composition factors, with multiplicity.
\end{propn}

\begin{pf}
First consider the Specht modules $S^\la$, $S^{\la'}$ for $\aaa n$, and consider the functor $M\mapsto M^\theta$ in \S\ref{duality}.  When $e=2$, it is quite easy to show that $(S^\la)^\theta\cong S^\la$ (and hence $(D^\la)^\theta\cong D^\la$ if $\la$ is $2$-regular).  Hence by the type $A$ version of Theorem \ref{duali} the Specht modules $S^\la$ and $S^{\la'}$ have the same composition factors.

Now consider Specht modules for $\bbb n$.  If $-Q$ is not a power of $q$, then by the previous paragraph and the results in \S\ref{morita} the Specht modules $S^{(\la,\mu)}$ and $S^{(\la',\mu)}$ have the same composition factors.  In particular, this applies when $Q$ is an indeterminate.  Now consider applying the second decomposition map in \S\ref{decompmap} to specialise $Q$: since a decomposition map is a function between Grothendieck groups which sends Specht modules to Specht modules, the result holds in general.
\end{pf}

\begin{pf}[Proof of Theorem \ref{main2}]
Suppose $\la\neq\varnothing\neq\mu$ and $S^{(\la,\mu)}$ is irreducible.  Then $(\la,\mu)$ is \infi t for every $t$ of the same parity as $r$.  By Propositions \ref{allirodd} and \ref{allireven}, this means that either $(\la,\mu)$ is \dable, or $r=0$ and $(\la,\mu)$ has the form $\left((2m),(1^{2n})\right)$ or $\left((1^{2n}),(2m)\right)$ for some $m,n>0$.  But by Proposition \ref{lacmu}, $S^{\left((2m),(1^{2n})\right)}$ has exactly the same composition factors as $S^{\left((1^{2m}),(1^{2n})\right)}$, and the latter Specht module is reducible, since the bipartition $\left((1^{2m}),(1^{2n})\right)$ is \infr{2m}.  A similar argument applies for the bipartition $\left((1^{2n}),(2m)\right)$.
\end{pf}

\end{document}